\documentclass[11pt]{amsart}

\usepackage{amsfonts,amsmath,amssymb,amsthm,mathrsfs,url,xcolor,latexsym,rawfonts,graphicx}

\usepackage{amsfonts}
\usepackage{pifont}
\usepackage{caption}
\usepackage{graphicx, subfig}
\usepackage{bm}
\usepackage{latexsym,amsmath,amssymb,cite,amsthm}
\usepackage{color,eucal,enumerate,mathrsfs}
\usepackage[normalem]{ulem}
\usepackage{array}
\usepackage{xy}
\usepackage{hyperref}
\usepackage[pagewise]{lineno}

\makeatletter
\@addtoreset{equation}{section}
\makeatother

\theoremstyle{plain}

\newtheorem{lemma}{Lemma}[section] 
\newtheorem{theorem}[lemma]{Theorem}

\newtheorem{remark}[lemma]{Remark} 

\newcommand{\Vol}{{\mbox{Vol}}}

\renewcommand{\det}{\mbox{det}}

\newcommand{\dist}{\mbox{dist}}

\usepackage{hyperref}

\usepackage{xcolor}

\newcommand{\beq}{\begin{equation}}
\newcommand{\eeq}{\end{equation}}
\newcommand{\beqs}{\begin{eqnarray*}}
\newcommand{\eeqs}{\end{eqnarray*}}
\newcommand{\beqn}{\begin{eqnarray}}
\newcommand{\eeqn}{\end{eqnarray}}
\newcommand{\beqa}{\begin{array}}
\newcommand{\eeqa}{\end{array}}

\begin{document}

	\title[The $L_p$ chord Minkowski problem for super-critical exponents]
	{The $L_p$ chord Minkowski problem for super-critical exponent}
	
	\author{Shibing Chen}
	\address{Shibing Chen, School of Mathematical Sciences,
		University of Science and Technology of China,
		Hefei, 230026, P.R. China.}
	\email{chenshib@ustc.edu.cn}

	\author[Q. Li]{Qi-Rui Li}
	\address{Qi-Rui Li, School of Mathematical Sciences, 
		Zhejiang University, Hangzhou 310027, China}
	\email{qi-rui.li@zju.edu.cn}

	\author[Y. Li]{Yuanyuan Li}
	\address{Yuanyuan Li, Institute for Theoretical Sciences, Westlake University, Hangzhou, 310030, China.}
	\email{lyyuan@westlake.edu.cn}

\thanks{
Research of the second author was supported by National Key R\&D Program of China (No.2022YFA1005500)
and Zhejiang Provincial NSFC (No.LR23A010002).
}	
	
\begin{abstract}
The $L_p$ chord Minkowski problem was recently introduced by Lutwak, Xi, Yang and Zhang, which seeks to determine the necessary and sufficient conditions for a given finite Borel measure  such that it is the $L_p$ chord measure of a convex body. In this paper, we solve the $L_p$ chord Minkowski problem for the super-critical exponents by combining a nonlocal Gauss curvature flow introduced in \cite{HHLW exi} and a topological argument developed in \cite{GLW2022}. Notably, we provide a simplified argument for the topological part.	
\end{abstract}

\maketitle

	\baselineskip=16.4pt
	\parskip=3pt

	\section{introduction}
	
	Recently, a new family of geometric measures were introduced by Lutwak, Xi, Yang and Zhang\cite{LXZY2022} by studying of a variational formula regarding intergral geometric invariants of convex bodies called chord integrals. Let $K \in \mathcal{K}^n,$ where $\mathcal{K}^n$ denotes the set of all convex bodies in $\mathbb{R}^n,$ the $q$-th chord integral $I_q(K)$ is defined by
	\begin{equation}\label{integral}
		I_q(K)=\int_{\mathcal{L}^n}|K\cap\ell|^q \mathrm{d}\ell,
	\end{equation}
	where $\mathcal{L}^n$ denotes the Grassmannian of 1-dimensional affine subspace of $\mathbb{R}^n,$ $|K\cap \ell|$ denotes the length of the chord $K\cap \ell,$ and the integration is with respect to Haar measure on the affine Grassmannian $\mathcal{L}^n,$ which is normalized to be a probability measure when restricted to rotations and to be $(n-1)$-dimensional Lebesgue measure when restricted to parallel translations.
	$$
	I_1(K)=V(K),\quad  I_0(K)=\frac{\omega_{n-1}}{n\omega_n}S(K), \quad I_{n+1}(K)=\frac{n+1}{\omega_n}V(K)^2,
	$$
	where $\omega_n$ denotes the volume of $n$-dimensional unit ball,
	and $V(K)$ denotes the volume of $K.$ One can see from the above fomula that the chord integrals include the convex body's volume and surface area as two special cases. These are Crofton’s volume formula, Cauchy’s integral formula for surface area, and the Poincar\'e-Hadwiger formula, respectively (see \cite{D. Ren, Santalo}).
	
	The chord measures and the Minkowski problems associated with chord measures were introduced in \cite{LXZY2022}. 
	They showed that the chord measures are the differentials of chord integrals and 
	solved the chord Minkowski problem except for the critical case of the Christoffel-Minkowski problem. 
Denote by $\mathcal{K}_o^n$ the set of all convex bodies containing the origin in the interior.
	For $K\in \mathcal{K}^n_{o}$ and $p,q\in \mathbb{R},$ the $L_p$ chord measures are defined by
	\begin{equation}\label{measure2}
		F_{p,q}(K,\eta)=\frac{2q}{\omega_n}\int_{\nu^{-1}_K(\eta)}(z\cdot\nu_{K}(z))^{1-p}\widetilde{V}_{q-1}(K,z)\mathrm{d} \mathcal{H}^{n-1}(z),\,\forall\text{  Borel set }\eta\subset\mathbb{S}^{n-1},
	\end{equation}
	where $\mathbb S^{n-1}$ is the unit sphere in $\mathbb R^n$,
	and $\widetilde{V}_{q-1}(K,z)$ is the ($q$-1)-th dual quermassintegral with respect to $z$. See \eqref{dual}.
	
When $q=1,$ $F_{p,1}(K,\cdot)$ corresponds the  $L_p$ surface area measure.  
The problem of characterizing the $L_p$ surface area measure is known as the $L_p$ Minkowski problem, 
which was first formulated and studied by Lutwak in \cite{L1}. 
Since then, the $L_p$ Minkowski problem with sub-critical exponent $p>-n$ has been extensively investigated,
see e.g. \cite{BLYZ13, BLYZ12(5), C-W2000}.
The case with super-critical exponent $p<-n$ was not resolved until recent work \cite{GLW2022},
where the authors introduced a topological method based on the calculation of the homology of a topological space of ellipsoids . 
For the classical Brunn-Minkowski theory and its recent developments, 
readers are referred to Schneider's monograph \cite{Schneider} and references therein.
	
The $L_p$ chord Minkowski problem posed by Lutwak, Xi, Yang and Zhang \cite{LXZY2022} is a problem of prescribing the $L_p$ chord measures. 
Given a finite Borel measure $\mu$ on $\mathbb{S}^{n-1},$ the $L_p$ chord Minkowski problem asks what are the necessary and sufficient conditions for $\mu$ such that it is the $L_p$ chord measure of a convex body $K\in \mathcal{K}^n_o$, namely
	\begin{equation}\label{original eq}
		F_{p,q}(K,\cdot)=\mu.
	\end{equation}
When $p=1,$ it is the chord Minkowski problem. When $q=1,$ it is the $L_p$ Minkowski problem.
When $\mu$ has a density $f$ that is an integrable nonnegtive function on $\mathbb{S}^{n-1}$,
the $L_p$ chord Minkowski problem is equivalent to solving the following Monge-Amp\`{e}re type equation
	\begin{equation}\label{MAeq}
		\mbox{det}(\nabla^2h+hI) =\frac{h^{p-1}f}{\widetilde{V}_{q-1}([h],\overline{\nabla} h)} \quad\text{on}\ \mathbb{S}^{n-1},
	\end{equation}
	where $h:\mathbb{S}^{n-1}\rightarrow \mathbb{R}$ is the support function of $K,$ $\nabla^2h$ is the covariant differentiation of $h$ with respect to an orthonormal frame on $\mathbb{S}^{n-1},$ $I$ is the unit matrix, 
$\overline{\nabla} h(x)=\nabla h(x)+h(x)x$ is the Euclidean gradient of $h$ in $\mathbb{R}^n$, 
and $\widetilde{V}_{q-1}([h],\overline{\nabla} h)$ is the ($q$-1)-th dual quermassintegral of the Wulff-shape $[h]$ 
with respect to the point $\overline{\nabla} h$.
	For detailed definitions, we refer readers to Section 2.
	
In \cite{LXZY2022}, Lutwak, Xi, Yang and Zhang found a sufficient condition for the origin-symmetric chord log-Minkowski problem
by studying the delicate concentration properties of cone-chord measures. 
Shortly afterward, Xi, Yang, Zhang and Zhao \cite{XYZZ2022} resolved the $L_p$ chord Minkowski problem
for $ p>1$ or $0<p<1$ under the origin-symmetric condition.
More recently, Guo, Xi and Zhao addressed the $L_p$ chord Minkowski problem for $0 \leq p < 1$ 
without any symmetry assumptions \cite{GXZ2023}. 
Subsequently, Li \cite{LYY exi} solved \eqref{MAeq} for $-n<p<0$ and $1\le q<n+1$,
and also provided a proof for the discrete $L_p$ chord Minkowski problem under the constraints $p<0$ and $q>0$. 
By a parabolic flow approach, Hu, Huang, Lu and Wang \cite{HHLW exi} obtained the existence of solutions to \eqref{MAeq}
when $f$ is positive, even and smooth, $p>-n$ and $p\neq 0$.
	
In this paper, we study \eqref{MAeq} for the super-critical exponents by using the method introduced in \cite{GLW2022}. 
As in \cite{HHL reg}, we study a Gauss curvature type flow
	\begin{equation}\label{flow1}
		\frac{\partial X}{\partial t}=-\frac{\omega_n f(\nu)\mathcal{K}\left\langle X,\nu\right\rangle ^p }{2q \widetilde{V}_{q-1}\left(\Omega_t,\overline{\nabla}(X\cdot \nu)\right)}\nu+X,
	\end{equation}
	with initial hypersurface $X(\cdot,0)=X_0(\cdot).$ 
Here $\mathcal{K}(\cdot, t)$ is the Gauss curvature of the convex hypersurface $\mathcal{M}_t,$ parametrized by smooth map 
$X(\cdot, t):\mathbb{S}^{n-1}\rightarrow \mathbb{R}^n,$ 
$\Omega_t=Cl(\mathcal{M}_t)$ is the convex body enclosed by $\mathcal{M}_t$,
and $\nu (\cdot, t)$ is the unit outer normal at $X(\cdot, t).$
Let $h(\cdot,t)$ be the support function of $\Omega_t$. 
Since the Gauss curvature of $\mathcal M_t$ is given by
\[
\mathcal{K} = \frac{1}{\det(\nabla^2h+hI)},
\]
it follows that
\begin{equation}\label{flow2}
 \partial_t h(x,t)=-\frac{\omega_n f(x)h(x,t)^p}{2q \widetilde{V}_{q-1}\left(\Omega_t,\overline{\nabla}h(x,t)\right)}\frac{1}{\det(\nabla^2h+hI)}
 +h(x,t), \ x\in\mathbb S^{n-1}.
\end{equation}
	
The dual quermassintegral $\widetilde{V}_{q-1}(K,z)$ is a nonlocal term and is difficult to deal with. 
Note that the ($q$-1)-th dual quermassintegral $\widetilde{V}_{q-1}(K, z)$ of $K$ with respect to $z\in \partial K$ is more delicate than the ($q$-1)-th dual quermassintegral $\widetilde{V}_{q-1}(K)$ of $K\in \mathcal{K}^n_o$. The main result of this paper is the following.
	
\begin{theorem}\label{A}
Let $p<-n-q+1$, $3< q<n+1$, and $\mu$ be a finite Borel measure on $\mathbb{S}^{n-1}$ with density $f$.
If $f\in C^{1,1}(\mathbb{S}^{n-1})$
and $\frac{1}{\Lambda}<f<\Lambda$ for some constant $\Lambda>0$,
then there exists a uniformly convex, positive, $C^{3,\alpha}$ solution to \eqref{MAeq}, where $\alpha \in (0,1).$
\end{theorem}	
	
	Applying an approximation argument, we can further obtain the existence of solutions when the density 
	$f \in L^{\infty}(\mathbb{S}^{n-1})$ satisfies $\frac{1}{\Lambda}<f<\Lambda$ for some $\Lambda>0.$
	
	\begin{theorem}\label{B}
		Let $p<-n-q+1$, $3<q<n+1,$ and $\mu$ be a finite Borel measure on $\mathbb{S}^{n-1}$ with density $f$.
		If $f \in L^{\infty}(\mathbb{S}^{n-1})$ and $\frac{1}{\Lambda}<f<\Lambda$ for some $\Lambda>0,$ 
		then there exists a strictly convex, positive, $C^{1,\alpha}$ weak solution to \eqref{MAeq}, where $\alpha \in (0,1).$
	\end{theorem}
	
Consider the following functional of convex bodies $\Omega \in \mathcal{K}_o$,
	\begin{equation}\label{functional}
		\mathcal{J}(\Omega)=I_q(K)-\frac{1}{p}\int_{\mathbb{S}^{n-1}}fh^p \mathrm{d}\sigma_{\mathbb{S}^{n-1}}.
	\end{equation}
We will show that \eqref{MAeq} is the Euler equation of this functional, and \eqref{flow1} constitutes a gradient flow associated with this functional. 
Hence if the flow \eqref{flow1} exists for all time and remains smooth and uniformly convex, 
then it deforms a initial hypersurface into a solution to \eqref{MAeq}.
The main difficulty of studying the flow \eqref{flow1} is the lack of uniform estimates for solutions.  
To address this challenge, we adopt a strategy akin to that employed in \cite{GLW2022}, albeit with simplifications in their proof for topological part.
	
To apply the method of \cite{GLW2022}, we first show that $\mathcal{J}(\Omega)$ is bigger than any given large constant in one of the following secnarios— the volume of $\Omega$ being sufficiently large or small, 
the eccentricity of $\Omega$ being sufficiently large, or the origin being close enough to the boundary of $\Omega.$  
Let $$\mathcal{A}_I:=\{E\in \overline{\mathcal{K}_o} \text{ is an ellipsoid in }\mathbb{R}^n: \omega_n\bar v\leq \Vol(E)\leq \omega_n\bar v^{-1} \text{ and }e_E\leq \bar{e}\},$$ where $\bar v, \bar{e}$ are appropriate constants.
Then, we construct a modified flow with initial data being an ellipsoid in $\mathcal{A}_I,$ similar to that in   \cite{GLW2022}. 
The key ingredient of \cite{GLW2022} is to show that there exists an initial data $\mathcal N$, which is an ellipsoid in $\mathcal{A}_I$, 
such that the flow \eqref{flow2} starting from $\mathcal N$ remains smooth and uniformly convex for all time $t\in [0,\infty)$.
For this end, a contradiction argument was employed.
Suppose such $\mathcal N$ does not exist. 
Then we have a retraction $\tilde{\Psi}$ from $\mathcal{A}_I$ to $\mathcal{P}$, the boundary of $\mathcal{A}_I$ given by 
\begin{equation}\label{pdef1}
			\mathcal{P} = \left\{ E \in \mathcal{A}_I: \text{ either } \operatorname{Vol}(E) =\omega_n \bar{v}, \text{ or } \operatorname{Vol}(E) = \frac{\omega_n}{\bar{v}}, \text{ or } e_E = \bar{e}, \text{ or } O \in \partial E \right\}.
\end{equation}
The original approach of \cite{GLW2022} then goes as follows. 
The existence of retraction $\tilde \Psi$ yields an injection from the homology group of $\mathcal{P}$ to that of $\mathcal{A}_I$. 
Therefore $\mathcal{P}$ possesses trivial homology since $\mathcal{A}_I$ is contractible, as shown in \cite[Lemmas 3.4 \& 3.5]{GLW2022}. 
The authors then calculated some homology group of $\mathcal{P}$ and showed that it is not trivial, 
see \cite[Proposition 3.6, Theorem 3.7]{GLW2022}.
A contradiction is thus arrived. 

The computation of the homology group is very delicate and involved. 
In this paper, we provide a simplified proof of this part by using the classical Brouwer fixed point theorem only, which might be helpful for readers. The key observations are as follows. First, since any ellipsoid \(E\) can be represented as \(E = A(B_1)\),
an affine transformation $A$ of the unit ball $B_1$.
We identify each affine transformation with an positive definite matrix.
Then \(\mathcal{A}_I\) is homeomorphic to \(\mathcal{E}_I \times B_1\), where
	\begin{equation*}
		\mathcal{E}_I = \{A \in M^{n \times n} \mid A \text{ is positive definite, } \bar{v} \leq \det(A) \leq \frac{1}{\bar{v}}, e_{A} \leq \bar{e}\},
	\end{equation*}
	for some \(0 < \bar{v} < 1\). 
Again we assume by contradiction that no nice initial data exists.
Then there is a retraction $\tilde{\Psi}$ from \(\mathcal{E}_I \times B_1\) to \(\mathcal{P}\). 
Denote
	\begin{equation*}
		\mathcal{D} := \{A \in M^{n \times n} \mid A \text{ is positive definite, } \|A\|_{\infty} \leq L, e_{A} \leq \bar{e}\},
	\end{equation*}
	where $L$ is a large constant such that $\mathcal{E}_I\subset \mathcal D.$
Our key observations are
\begin{itemize}
\item \(\mathcal{D}\) is convex;
\item there is a retraction \(\Phi\) from \(\mathcal{D}\times B_1\) to \(\mathcal{E}_I \times B_1\);
\item we can construct a mapping \(g\) from \(\mathcal{P}\) to itself such that \(g\) has no fixed points. 
\end{itemize}
Let \(i: \mathcal{P} \rightarrow \mathcal D\times B_1\) be the inclusion map. Then
	\begin{equation*}
		G = i \circ g \circ \tilde{\Psi} \circ \Phi: \mathcal D\times B_1 \rightarrow \mathcal D\times B_1
	\end{equation*}
is a continuous map without fixed points, contradicting the Brouwer fixed point theorem.

The paper is organized as follows. In Section 2, we present some basic concepts in the theory of convex bodies and integral geometry and 
recall some relevant theorems from the literature. In Section 3, we derive the $C^2$-estimates of solutions to the flow \eqref{flow2}
by assuming the $C^0 \& C^1$-estimates. 
In Section 4, we first prove the monotonicity and give some estimates of the functional \eqref{functional}, 
and then introduce a modified flow associated to \eqref{flow1}. 
Section 5 is dedicated to proving Theorems \ref{A} and \ref{B}.

	\section{Preliminaries}
	In this section, we introduce necessary notations and collect relevant results from the literatures that will be useful for the subsequent analysis. 
	
Let $x\cdot y$ be the inner product of $x,y\in \mathbb{R}^n,$ and $|x|=\sqrt{x\cdot x}$ be the Euclidean norm of $x.$ 
A convex body $K$ is a compact convex subset of $\mathbb{R}^n$ with non-empty interior. 
Denote by $\mathcal{K}^n$ the set of all convex bodies in $\mathbb{R}^n$,
and by $\mathcal{K}_o$ the set of convex bodies that contains the origin in the interior.
	For a continuous function $h : \mathbb{S}^{n-1}\rightarrow (0, \infty),$ the Wulff shape of $h$ is the convex body
	$$
	[h]:=\big\{x\in\mathbb{R}^n: x\cdot u\leq h(u)\text{ for all } u\in\mathbb{S}^{n-1}\big\}.
	$$
	
	Let $K\in \mathcal{K}^n,$ and $h_K(v):=\max\{v\cdot x,x\in K\},$ $\rho_K(u):=\max\{\lambda:\lambda u\in K\}$ are the support function and the radial function of convex body $K$ defined from $ \mathbb{S}^{n-1}\rightarrow \mathbb{R}.$ We write the support hyperplane of $K$ with the outer unit normal $v$ as
	$$
	H_K(v):=\left\{x\in \mathbb{R}^n: x\cdot v=h_K(v) \right\},
	$$
	and the half-space $H^{-}(K,v)$ in direction $v$ is defined by
	$$
	H^{-}_K(v):=\left\{x\in \mathbb{R}^n: x\cdot v\leq h_K(v) \right\}.
	$$  
	Denote $\partial K$ as the boundary of $K$, that is, $\partial K:=\{\rho_K(u)u:u\in \mathbb{S}^{n-1}\}.$ The spherical image $\nu_{K}:\partial K\rightarrow \mathbb{S}^{n-1}$ is given by
	\begin{equation}\label{Guass}
		\nu_K(\{x\}):=\{v\in\mathbb{S}^{n-1}:x\in H_K(v)\}.
	\end{equation}
	Let $\sigma_K\subset \partial K$ denote the set of all points $x\in \partial K,$ such that the set $\nu_K(\{x\})$ contains more than one element. From \cite[P. 84]{Schneider}, we have $\mathcal{H}^{n-1}(\sigma_K)=0$. The function
	$$
	\nu_K: \partial K \backslash \sigma_K \longrightarrow \mathbb{S}^{n-1},
	$$
	defined by letting $\nu_K(x)$ be the unique element in $\nu_K(\{x\})$ for each $x \in \partial K \backslash \sigma_K$, is called the spherical image map of $K$ and is known to be continuous \cite[Lemma 2.2.12]{Schneider}. 
	
	
	Let $K \in \mathcal{K}^n$. For $z \in \operatorname{int} K$ and $q \in \mathbb{R}$, the $q$th dual quermassintegral $\widetilde{V}_q(K, z)$ of $K$ with respect to $z$ is defined by
	\begin{equation}\label{dual}
		\widetilde{V}_q(K, z):=\frac{1}{n} \int_{\mathbb{S}^{n-1}} \rho_{K, z}(u)^q \mathrm{d} u,
	\end{equation}
	where $\rho_{K, z}(u):=\max \{\lambda>0: z+\lambda u \in K\}$ is the radial function of $K$ with respect to $z$ defined from $\mathbb{S}^{n-1}$ to $\mathbb{R}.$ When $z \in \partial K$, $\widetilde{V}_q(K, z)$ is defined in the way that the integral is only over those $u \in \mathbb{S}^{n-1}$ such that $\rho_{K, z}(u)>0$. In other words,
	$$
	\widetilde{V}_q(K, z):=\frac{1}{n} \int_{\{u\in \mathbb{S}^{n-1}:\rho_{K, z}(u)>0\}} \rho_{K, z}(u)^q \mathrm{d} u \ ,\, z \in \partial K .
	$$
	When $q>-1,$ for $\mathcal{H}^{n-1}$-almost all $z \in \partial K$, we have
	\begin{equation}\label{dual ray}
		\widetilde{V}_q(K, z)=\frac{1}{2 n} \int_{\mathbb{S}^{n-1}} X_K(z, u)^q \mathrm{d} u,
	\end{equation}
	where the parallel $X$-ray of $K$ is the nonnegative function on $\mathbb{R}^n \times \mathbb{S}^{n-1}$ defined by
	$$
	X_K(z, u)=|K \cap(z+\mathbb{R} u)|, \quad z \in \mathbb{R}^n, \, u \in \mathbb{S}^{n-1} .
	$$
	When $q>0$, the dual quermassintegral is the Riesz potential of the characteristic function, that is,
	$$
	\widetilde{V}_q(K, z)=\frac{q}{n} \int_K|x-z|^{q-n} \mathrm{d} x.
	$$
	Note that this immediately allows for an extension of $\widetilde{V}_q(K, \cdot)$ to $\mathbb{R}^n$. An equivalent definition via radial function can be found in \cite{LXZY2022}. By a change of variables, we obtain:
	$$
	\widetilde{V}_q(K, z)=\frac{q}{n} \int_{K-z}|y|^{q-n} \mathrm{d} y,
	$$
	where $K-z:=\{x\in \mathbb{R}^n:x=y-z\text{ for some }y\in K\}.$ Indeed, when $q>0$, the integrand $|y|^{q-n}$ being locally integrable, it can be inferred that the dual quermassintegral $\widetilde{V}_q(K, z)$ is continuous in $z$. Let $K \in \mathcal{K}^n$. 
	When $z \in \partial K$, then either $\rho_{K, z}(u)=0$ or $\rho_{K, z}(-u)=0$ for almost all $u \in \mathbb{S}^{n-1}$, and thus
	\begin{equation}\label{ray and radial}
		X_K(z, u)=\rho_{K, z}(u), \quad \text { or } X_K(z, u)=\rho_{K, z}(-u), \quad z \in \partial K,
	\end{equation}
	for almost all $u \in \mathbb{S}^{n-1}$. 
	
	As presented before, let $q>-1$ and $K\in \mathcal{K}^n,$ the $q$-th chord integral of $K$ is given by
	$$
	I_q(K)=\int_{\mathcal{L}^n}|K\cap\ell|^q \mathrm{d}\ell,
	$$
	where $\mathcal{L}^n$ denotes the Grassmannian of 1-dimensional affine subspace of $\mathbb{R}^n,$ $|K\cap \ell|$ denotes the length of the chord $K\cap \ell,$ and the integration is with respect to Haar measure on the affine Grassmannian $\mathcal{L}^n.$ For $q>0,$ the chord integral can be written as the integral of dual quermassintegrals in $z\in K:$
	$$
	I_q(K)=\frac{q}{\omega_n}\int_{K}\widetilde{V}_{q-1}(K, z) \mathrm{d}z.
	$$
	When $q \geq 0 $, the chord integral $I_q(K)$ can be represented as follows:
	$$
	I_q(K)=\frac{1}{n \omega_n} \int_{\mathbb{S}^{n-1}} \int_{K|u^{\bot}} X_K(x, u)^q \mathrm{d} x \mathrm{d} u.
	$$
	When $q>1,$ the chord integral can be recognized as Riesz potential:
	$$
	I_q(K)=\frac{q(q-1)}{n\omega_n}\int_{K}\int_K|x-z|^{q-n-1} \mathrm{d} x\mathrm{d}z.
	$$
	An elementary property of the functional $I_q$ is its homogeneity. If $K \in \mathcal{K}^n$ and $q >-1$, then
	\begin{equation}\label{chord homogen}
		I_q(t K)=t^{n+q-1} I_q(K),
	\end{equation}
	for $t>0$. By compactness of $K$, it is easy to see that the chord integral $I_q(K)$ is finite whenever $q \geq 0$.
	Let $K \in \mathcal{K}^n$ and $q>0$, the chord measure $F_q(K, \cdot)$ is a finite Borel measure on $\mathbb{S}^{n-1}$, which can be expressed as:
	\begin{equation}\label{chord measure}
		F_q(K, \eta)=\frac{2 q}{\omega_n} \int_{\nu_K^{-1}(\eta)} \widetilde{V}_{q-1}(K, z) \mathrm{d} \mathcal{H}^{n-1}(z), \quad \text { for each Borel } \eta \subset \mathbb{S}^{n-1}.
	\end{equation}
	The mapping $\nu_{K}$ of $K$ is almost everywhere defined on $\partial K$ with respect to the $(n-1)$-dimensional Hausdorff measure, owing to the convexity of $K$. The chord measure $F_q(K, \cdot)$ is significant as it is obtained by differentiating the chord integral $I_q$ in a certain sense, as shown in \eqref{diff}. Chord measures inherit its translation invariance and homogeneity (of degree $n+q-2$) from chord integrals. And it is evident that the chord measure $F_q(K, \cdot)$ is absolutely continuous with respect to the surface area measure $S_{n-1}(K, \cdot)$. In \cite[Theorem 4.3]{LXZY2022}, it was demonstrated that:
	\begin{equation}\label{chord integral}
		I_q(K)=\frac{1}{n+q-1} \int_{\mathbb{S}^{n-1}} h_K(v) \mathrm{d} F_q(K, v).
	\end{equation}
	
	For each $p \in \mathbb{R}$ and $K \in \mathcal{K}_o^n$, the $L_p$ chord measure $F_{p, q}(K, \cdot)$ is defined as follows:
	\begin{equation}\label{lp chord}
		\mathrm{d} F_{p, q}(K, v)=h_K(v)^{1-p} \mathrm{d} F_q(K, v).
	\end{equation}
	We also have an important property of $F_{p,q},$ its homogeneity, namely
	\begin{equation}\label{lp chord homo}
		F_{p,q}(tK,\cdot)=t^{n+q-p-1}F_{p,q}(K,\cdot)
	\end{equation}
	for each $t>0.$
	
	From Theorem 2.2 in \cite{XYZZ2022}, we know that if $K_i\in \mathcal{K}^{n}_o \rightarrow K_0\in \mathcal{K}^{n}_o,$ then the chord measure $F_q(K_i,\cdot)$ converges to $F_q(K,\cdot)$ weakly. Hence, one can immediately obtain that 
	\begin{equation}\label{weak converg}
		F_{p,q}(K_i,\cdot) \rightarrow F_{p,q}(K,\cdot) \text{ weakly. }
	\end{equation}
	It was shown in \cite{LXZY2022} that the differential of the chord integral $I_q$ with respect to the $L_p$ Minkowski combinations leads to the $L_p$ chord measure: for $p \neq 0$,
	$$
	\left.\frac{\mathrm{d}}{\mathrm{d} t}\right|_{t=0} I_q\left(K+_p t \cdot L\right)=\frac{1}{p} \int_{\mathbb{S}^{n-1}} h_L^p(v) \mathrm{d} F_{p, q}(K, v),
	$$
	where $K+_p t \cdot L$ is the $L_p$ Minkowski combination between $K$ and $L.$
	
	To prove the monotonicity of the functional \eqref{functional}, we need the following variational formula for chord integral. 
	\begin{theorem}[Theorem 5.5 in \cite{LXZY2022}]
		Let $q>0,$ and $\Omega$ be a compact subset of $\mathbb{S}^{n-1}$ that is not contained in any closed hemisphere. Suppose that $g:\Omega\rightarrow \mathbb{R}$ is continuous and $h_t: \Omega \rightarrow (0,\infty)$ is a family of continuous functions given by
		$$
		h_t=h_0+tg+o(t,\cdot),
		$$
		for each $t\in(-\delta,\delta)$ for some $\delta>0.$ Here $o(t,\cdot)\in C(\Omega)$ and $o(t,\cdot)/t$ tends to $0$ uniformly on $\Omega$ as $t\rightarrow 0.$ Let $K_t$ be the Wulff shape generated by $h_t$ and $K$ be the Wulff shape generated by $h_0.$ Then,
		\begin{equation}\label{diff}
			\left.\frac{\mathrm{d}}{\mathrm{d} t}\right|_{t=0}I_q(K_t)=\int_{\Omega}g(v)\mathrm{d}F_q(K,v).
		\end{equation}
	\end{theorem}

	In the subsequent analysis, we would frequently utilize the lower bound of $I_q.$  
	\begin{lemma}\label{elliptic estimate}
		Suppose $q>1,$ if $E$ is an ellipsoid in $\mathbb{R}^n$ given by
		$$
		E=\left\{x \in \mathbb{R}^n: \frac{ x_1^2}{a_1^2}+\cdots+\frac{x_n^2}{a_n^2} \leq 1\right\}
		$$
		with $0<a_1 \leq a_2 \leq \cdots \leq a_n$, then we have
		\begin{equation}\label{Iq below}
			I_q(E) \geq c_n a_2 \cdots a_n a_1^q
		\end{equation}
		for some positive constant $c_n$ depending only on $n.$
	\end{lemma}
	\begin{proof}
		Since
		$$
		I_q(E)=\frac{1}{n \omega_n} \int_{S^{n-1}} \int_{E|u^{\bot}} X_{E}(x, u)^q \mathrm{~d} x \mathrm{~d} u, \quad q \geq 0 .
		$$
		where $E|u^{\bot}$ denotes the projection of $E$ onto $u^{\bot}.$ For $q>1,$ we have 
		\begin{equation}\label{jess1}
			I_q(E)\geq \frac{1}{n \omega_n} \int_{S^{n-1}} V(E)^q V_{n-1}(E|u^{\bot})^{1-q}\mathrm{d}u
		\end{equation}
		Indeed, Jessen inequality gives
		\begin{small}
			$$
			\frac{1}{V_{n-1}(E|u^{\bot})}\int_{E_i|u^{\bot}} X_{E}(x, u)^q \mathrm{~d} x\geq \left(\frac{1}{V_{n-1}(E|u^{\bot})}\int_{E|u^{\bot}} X_{E}(x, u) \mathrm{~d} x\right)^q=\left(\frac{V(E)}{V_{n-1}(E|u^{\bot})}\right)^q.
			$$
		\end{small}
		Recall that $V(E)=\omega_{n} \prod_{1}^{n}a_i,$ and it is straightforward to check that $V_{n-1}(E|u^{\bot})\leq C_n\frac{V(E)}{a_1}\leq C_n\omega_{n}\prod_{2}^{n}a_i$ for some constant $C_n$ depending only on $n.$
		Hence by \eqref{jess1} we have
		$$
		I_q(E) \geq c_n a_2 \cdots a_n a_1^q
		$$
		some positive constant $c_n$ depending only on $n.$
	\end{proof}
	

\section{A priori estimates for solutions to the Gauss curvature flow}
	In this section, we establish the derivative estimates for solutions to \eqref{flow2}. 
	
\begin{theorem}\label{C2}
Let $f$ be a positive and $C^{1,1}$-smooth function on $\mathbb{S}^{n-1}$, $p<-n-q+1$ and $3< q<n+1$. 
Let $h(\cdot, t)$ be a positive, smooth and uniformly convex solution to \eqref{flow2} for $t \in [0, T)$. 
Assume that
		\begin{equation}\label{C0C1}
			\begin{aligned}
				1 / C_0 \leq h(x, t) &\leq C_0, \\
				|\nabla h|(x, t) &\leq C_0,
			\end{aligned}
		\end{equation}
		for all $(x, t) \in \mathbb{S}^{n-1} \times[0, T).$ Then
		\begin{equation}\label{2th derivative}
			C^{-1} I \leq\left(\nabla^2 h+h I\right)(x, t) \leq C I \quad \forall(x, t) \in \mathbb{S}^{n-1} \times[0, T),
		\end{equation}
for some constant $C>0$ depending only on $n, p, q,C_0, \min _{\mathbb{S}^{n-1}} f,\|f\|_{C^{1,1}\left(\mathbb{S}^{n-1}\right)}$, 
and the initial condition $h(\cdot, 0)$.
	\end{theorem}

By approximation, we may assume directly that $f$ is $C^2$-smooth. 
The proof of Theorem \ref{C2} uses similar ideas as in \cite{GLW2022,HHL reg}.

	\begin{proof}[\bfseries{Proof of Theorem \ref{C2}}]
Let $\mathcal M_t$ be the boundary of the Wulff shape $[h(\cdot,t)]$. Then $\mathcal M_t$ is evolved by \eqref{flow1}.
The proof is divided into two steps. 

\noindent		
{\bfseries{Step 1}}: 
$\max_{\mathbb S^{n-1}\times[0,T)}\frac{1}{\det (\nabla^2 h+hI)}\leq C$.
		
Recall that the principal radii of curvature of $\mathcal M_t$ are eigenvalues of the matrix 
		$$
		b_{ij}=h_{ij}+h \delta_{ij},
		$$
and so the Gauss curvature $\mathcal{K}$ of $\mathcal{M}_t$ is
\[
\mathcal{K} = \frac{1}{\det b_{ij}}  = \frac{1}{\det (\nabla^2h+hI)} .
\]
Consider the auxiliary function:
\begin{equation}\label{def Q}
		Q=-\frac{h_t}{h-\varepsilon_0}=\frac{1}{h-\varepsilon_0}\left(\frac{\omega_n f \mathcal{K} h ^p }{2q \widetilde{V}_{q-1}\left(\Omega_t,\overline{\nabla}h\right)}-h\right),
\end{equation}
		where $\varepsilon_0=\frac{1}{2} \min _{\mathbb{S}^{n-1} \times[0, T)} h>0$ 
		and $\Omega_t$ is the convex body given by $h(\cdot,t).$ 
		By \eqref{C0C1},
		\begin{equation}\label{Vq estimate}
			\frac1C\leq \widetilde{V}_{q-1}\left(\Omega_t,\overline{\nabla}h\right)\leq C,
		\end{equation}
	and
	\[
	\frac1CQ\leq \mathcal{K}\leq C(Q+1).
	\]
In the sequel, we always use $C$ to denote a positive constant
which depends only on $n,p,q,C_0$, $\min_{\mathbb S^{n-1}} f$, $\|f\|_{C^{1,1}(\mathbb S^{n-1})}$ and $h(\cdot,0)$, 
but it may change from line to line.
To complete Step 1, it suffices to estimate $Q$ from above.
				
		For any given $T'\in(0,T)$, let us assume
\[
Q(x_0,t_0) = \max_{\mathbb{S}^{n-1} \times\left[0, T^{\prime}\right]}Q.
\]
If $t_0=0$, then $\max _{\mathbb{S}^{n-1} \times\left[0, T^{\prime}\right]} Q=\max _{\mathbb{S}^{n-1}} Q(\cdot, 0)$ and we are done. 
Suppose $t_0>0$. Then
		$$
		0=\nabla_i Q|_{(x_0,t_0)}=-\frac{h_{t i}}{h-\varepsilon_0}+\frac{h_t h_i}{\left(h-\varepsilon_0\right)^2},
		$$
which gives, at $(x_0,t_0)$,
		\begin{equation}\label{Q C1}
			h_{t i}=-Q h_i=\frac{h_t h_i}{h-\varepsilon_0}.
		\end{equation}
We also have
		\begin{equation*}
			\begin{aligned}
				0 \geq \nabla_{i j} Q|_{(x_0,t_0)}& =\frac{-h_{t i j}}{h-\varepsilon_0}+\frac{h_{t i} h_j+h_{tj}h_i+h_t h_{i j}}{\left(h-\varepsilon_0\right)^2}-\frac{2 h_t h_i h_j}{\left(h-\varepsilon_0\right)^3} \\
				& =\frac{-h_{t i j}}{h-\varepsilon_0}+\frac{h_t h_{i j}}{\left(h-\varepsilon_0\right)^2}\\
				& = \frac{-h_{t i j}-Qh_{i j}}{h-\varepsilon_0},
			\end{aligned}
		\end{equation*}
which yields, at $(x_0,t_0)$,
\begin{equation}\label{iiQ}
 h_{tij}\ge -Q h_{ij}.
\end{equation}

Differentiating \eqref{def Q} with respect to $t$ gives
		\begin{equation}\label{ptQ}
			\begin{aligned}
				0\leq \partial_t Q|_{(x_0,t_0)}&= \frac{-h_{t t}}{h-\varepsilon_0}+\frac{h_t^2}{\left(h-\varepsilon_0\right)^2} \\
				&=\frac{\omega_n f}{2q(h-\varepsilon_0)}\left[ \frac{h^p}{\widetilde{V}_{q-1}\left(\Omega_t, \bar{\nabla} h\right)}\partial_t \mathcal{K}+\frac{\mathcal{K}}{\widetilde{V}_{q-1}\left(\Omega_t, \overline{\nabla} h\right)}\partial_t (h^p)\right. \\
				&\left.+\mathcal{K} h^p \frac{d \widetilde{V}_{q-1}^{-1}\left(\Omega_t, \overline{\nabla} h\right)}{d t}\right]+Q+Q^2 .
			\end{aligned}
		\end{equation}
We next estimate the terms of \eqref{ptQ}.
		Let $\{b^{i j}\}$ be the inverse matrix of $\{b_{i j}\}.$ By a rotation of coordinates, 
		we may assume that $\{b_{i j}\}$ is diagonal at $(x_0,t_0)$. Then
		\begin{equation}\label{eq of bij}
			\sum b^{i i}\geq (n-1)\left(\prod b^{i i}\right)^{\frac{1}{n-1}}=(n-1) \mathcal{K}^{\frac{1}{n-1}}\ge \frac1CQ^{\frac{1}{n-1}}.
		\end{equation}
		By \eqref{iiQ} and \eqref{eq of bij}, we have that
		\begin{equation}\label{ptK}
			\begin{aligned}
				\partial_t \mathcal{K}|_{(x_0,t_0)} 
				&= -\mathcal{K} b^{i j}(h_{tij}+h_t \delta_{i j})\\
				&\le \mathcal{K} b^{i j}( Q (b_{ij}-h\delta_{ij}) -h_t \delta_{i j})\\
				& = \mathcal{K} Q \left(n-1-\varepsilon_0 \sum b^{i i}\right) \\
				& \leq -\varepsilon_0CQ^{2+\frac{1}{n-1}} +C(Q^2+1).
			\end{aligned}
		\end{equation}
Also
		\begin{equation}\label{pthp}
			\partial_t (h^p)=p h^{p-1} \partial_t h\leq CQ.
		\end{equation}
Denote $z=z(x,t) = \overline\nabla h(x,t)$ and $z_0 = \overline\nabla h(x_0,t_0)$.
Let
\[
 h_{z}(\xi,t) =  h(\xi,t)-\xi\cdot z(x,t), \ \text{for} \ \xi\in\mathbb S^{n-1}.
\]
Let $S^+_{z} = \{u\in\mathbb S^{n-1}: \rho_{z}(u,t)>0\}$, where $\rho_z(u,t) = \max\{\lambda>0:z+\lambda u\in \Omega_t\}$.
For every $u\in S^+_{z} $, one sees 
\begin{equation}\label{thm3.1:tt1}
\langle \xi(u,z), u\rangle \rho_{z} (u,t)= h_{z}(\xi(u,z),t),
\end{equation}
where $\xi(u,z)$ is the unit outer normal of $\mathcal M_t$ at $z(x,t)+u\rho_{z} (u,t)$.
Denote $\dot\xi (u,z)= \frac{d}{dt}\xi(u,z(x,t))$. We then differentiate \eqref{thm3.1:tt1} and find for every $u\in S^+_{z} $,
\begin{eqnarray*}
\frac{\frac{d}{dt}\rho_{z} (u,t)}{\rho_{z} (u,t)}
&=&\frac{\frac{d}{dt}h_z(\xi(u,z),t)}{h_z(\xi(u,z),t)}-\frac{\langle  \dot\xi(u,z), u\rangle}{\langle \xi(u,z), u\rangle}\\
&=&\frac{\partial_t h(\xi(u,z),t) -\xi(u,z)\cdot \partial_tz(x,t) +\langle\nabla_{\xi} h_z, \dot\xi(u,z)\rangle}{h_z(\xi(u,z),t)}
-\frac{\langle \dot\xi(u,z), u\rangle}{\langle \xi(u,z), u\rangle}\\
&=&\frac{\partial_t h(\xi(u,z),t) -\xi(u,z)\cdot \partial_tz(x,t) }{h_z(\xi(u,z),t)}
+\frac{\langle\nabla_{\xi} h_z-u \rho_{z} (u,t), \dot\xi(u,z)\rangle}{h_z(\xi(u,z),t)}\\
&=&\frac{\partial_t h(\xi(u,z),t) -\xi(u,z)\cdot \partial_tz(x,t) }{h_z(\xi(u,z),t)}
+\frac{\langle-h_z(\xi(u,z),t)\xi(u,z), \dot\xi(u,z)\rangle}{h_z(\xi(u,z),t)} \\
&=&\frac{\partial_t h(\xi(u,z),t) -\xi(u,z)\cdot \partial_tz(x,t) }{h_z(\xi(u,z),t)}.
\end{eqnarray*}
The last equality uses $\xi(u,z)\cdot \dot \xi(u,z) = 0$.
Hence, for every $u\in S^+_{z_0} $,
\begin{eqnarray}\label{rho diff}
 \frac{d}{dt}\rho_{z} (u,t_0)\big|_{z=z_0} 
 &=& \frac{\partial_t h(\xi(u,z_0),t_0) -\langle \xi(u,z_0), \nabla h_t+\partial_thx_0\rangle}{\xi(u,z_0)}\notag\\
 &=& \frac{\partial_t h(\xi(u,z_0),t_0) +Q\langle \xi(u,z_0), z_0-\varepsilon_0x_0\rangle}{\xi(u,z_0)\cdot u}
\end{eqnarray}		
where \eqref{Q C1} is used in the last equality.
By \eqref{dual}, \eqref{rho diff}, we obtain
\begin{eqnarray}\label{ptV}
\frac{d}{d t} \widetilde{V}_{q-1}^{-1}(\Omega_t,\overline\nabla h) \big|_{(x_0,t_0)}&=& 
-\frac{q-1}{n \widetilde V^2_{q-1}(\Omega_t,\overline\nabla h)} \int_{S^+_{z_0}} \rho_z^{q-2}\frac{d}{dt}\rho_{z}(u,t_0)\big|_{z=z_0}du\notag\\
&\le& C Q\int_{S^+_{z_0}} \frac{\rho_{z_0}^{q-2}(u,t_0)}{\xi(u,z_0)\cdot u} du.
\end{eqnarray}	
Inserting \eqref{ptK}, \eqref{pthp} and \eqref{ptV} into \eqref{ptQ}, we obtain 
\begin{eqnarray}\label{tQ2}
0\le-\varepsilon_0 CQ^{2+\frac{1}{n-1}}+C(Q^2+1)+CQ^2\int_{S^+_{z_0}} \frac{\rho_{z_0}^{q-2}(u,t_0)}{\xi(u,z_0)\cdot u} du.
\end{eqnarray}	

Let $v = v(u) \in\mathbb S^{n-1}$ be such that
\[
 \rho(v(u),t_0) v(u) = z_0+\rho_{z_0}(u,t_0)u.
\]
Recall that if $x= x(v)$ and $\xi = \xi(u)$ are respectively the unit outer normals of $\mathcal M_{t_0}$ at $ \rho(v,t_0) v$ and $\rho_{z_0}(u,t_0)u$,
then the following variable change formulas are known, see e.g. \cite{LSW Jems},
\[
\frac{dx}{dv}= \frac{\rho^n(v,t_0)}{h\det(\nabla^2h+hI)(x(v),t_0)} \quad\text{and}\quad 
\frac{d\xi} {du}= \frac{\rho_{z_0}^n(u,t_0)}{h_{z_0}\det(\nabla^2h_{z_0}+h_{z_0}I)(\xi(u),t_0)}.
\]
It follows that $v=v(u)$ satisfies the following variable change formula
\[
\frac{dv}{ du}= \frac{h(\xi(u),t_0)\rho_{z_0}^n(u,t_0)}{h_{z_0}(\xi(u),t_0)\rho^n(v(u),t_0)}.
\]
As a result,
\begin{eqnarray*}
\int_{S^+_{z_0}} \frac{\rho_{z_0}^{q-2}(u,t_0)}{\xi(u,z_0)\cdot u} du 
&=& \int_{u(v)\in  S^+_{z_0}} \rho_{z_0}^{q-1-n}(u(v),t_0)\frac{\rho^n(v(u),t_0)}{h(\xi(u),t_0)}dv\notag\\
&\le &C \int_{u(v)\in S^+_{z_0}} \rho_{z_0}^{q-1-n}(u(v),t_0)dv\notag\\
&\le&C \int_{\mathbb S^{n-1}} |\rho(v,t_0)v-z_0|^{q-1-n}dv\notag\\
&\le &C,
\end{eqnarray*}
where we use $q>2$ in the last inequality. Plugging this in \eqref{tQ2}, we get at $(x_0,t_0)$,
\[
 0\le -\varepsilon_0 CQ^{2+\frac{1}{n-1}}+C(Q^2+1).
\]
Therefore $\max_{\mathbb S^{n-1}\times [0,T']}Q \leq C$. Since $C$ is independent of $T'$, we are through.

\bigskip
\noindent	
{\bfseries{Step 2}}: $\nabla^2 h(x,t)+h(x,t)I\le CI$ for all $(x,t)\in\mathbb S^{n-1}\times[0,T)$.

Consider the following auxiliary function
		\begin{equation}\label{E x,t}
			E(x, t)=\log \lambda_{\max }\left(\left\{b_{i j}(x,t)\right\}\right)-A \log h(x, t)+B|\nabla h(x,t)|^2,
		\end{equation}
		where $A$ and $B$ are positive constants to be determined later, and $\lambda_{\max }\left(\{b_{ij}\}\right)$ denotes the maximal eigenvalue of $\{b_{ij}\}$. For any given $0<T^{\prime}<T$, suppose
\[
\max_{\mathbb S^{n-1}\times[0,T']}E=E\left(x_0, t_0\right)
\] 
Again, we assume w.l.o.g. $t_0>0.$ 
By a rotation of coordinates, we also assume that $\{b_{ij}\}\left(x_0, t_0\right)$ is diagonal, 
and $\lambda_{\max }\left(\{b_{ij}(x_0, t_0)\}\right)=b_{11}(x_0, t_0)$. 
Hence we only need to derive the upper bound for the following quantity
$$
		E(x, t)=\log b_{11}-A\log h(x, t)+B|\nabla h(x,t)|^2.
$$
At $(x_0,t_0)$, we have
\begin{equation}\label{T3.1:cc1}
0=\nabla_iE = b^{11} b_{11i} - A\frac{h_i}{h}+2Bh_i h_{ii},
\end{equation}
where $(b^{ij})$ is the inverse of $(b_{ij})$; and
\begin{equation}\label{T3.1:cc2}
0\ge\nabla_{ii}E = b^{11}\nabla_{ii}^2b_{11}-(b^{11})^2(\nabla_ib_{11})^2 
-A\Big(\frac{h_{ii}}{h}-\frac{h_i^2}{h^2}\Big) + 2B\Big(h_{ii}^2+\sum_kh_k h_{kii}\Big).
\end{equation}
By \eqref{flow2}, we have
\begin{equation}\label{log func}
\log \left(h-h_t\right) =-\log \det(\nabla^2 h +hI)+\psi(x,t),
\end{equation}
where
\[
		\psi(x, t):=\log \left(\frac{\omega_n f h ^p }{2q \widetilde{V}_{q-1}\left(\Omega_t,\overline{\nabla}h\right)}\right).
\]
Since the equation \eqref{log func} is of the same type as in \cite[Lemma 5.2]{GLW2022}, 
we follow the same computation of that paper and obtain
at $(x_0,t_0)$ that, see (5.13) in \cite[Lemma 5.2]{GLW2022},
\begin{eqnarray}\label{T3.1:cc3}
\frac{\partial_t E}{h-h_t} &\le &b^{11}\Big[\sum b^{ii}(\nabla^2_{ii}b_{11} -b_{11}+b_{ii}) -\sum b^{ii}b^{jj}(\nabla_1 b_{ij})^2\Big]\notag\\
&&-b^{11}\nabla^2_{11}\psi+\frac{1-A}{h-h_t} + \frac Ah+2B\frac{\sum h_kh_{kt}}{h-h_t}.
\end{eqnarray}
Inserting \eqref{T3.1:cc2} into \eqref{T3.1:cc3}, we find at $(x_0,t_0)$ that
\begin{eqnarray*}
\frac{\partial_t E}{h-h_t}&\le& \sum b^{ii}\Big[ \frac{(\nabla_i b _{11})^2}{b^2_{11}} +A\Big(\frac{h_{ii}}{h} - \frac{h_i^2}{h}\Big)
-2B\Big(h_{ii}^2+\sum h_k h_{kii}\Big)\Big]\\
&&-b^{11}\sum b^{ii}b^{jj}(\nabla_1b_{ij})^2 - \frac{\nabla_{11}^2\psi}{b_{11}}
+\frac{1-A}{h-h_t}+\frac Ah+2B\frac{\sum h_k h_{kt}}{h-h_t}.
\end{eqnarray*}
Using $\sum b^{ii}b^{11}(\nabla_ib_{11})^2\le \sum b^{ii}b^{jj} (\nabla_1b_{ij})^2$, we further obtain
\beqn\label{ptE}
\frac{\partial_t E}{h-h_t}&\le& \sum b^{ii}A\Big(\frac{h_{ii}}{h} - \frac{h_i^2}{h}\Big)
-2B\sum b^{ii}h_{ii}^2+2B \sum h_k\Big(\sum -b^{ii} h_{kii}+\frac{h_{kt}}{h-h_t}\Big)\notag\\
&&- \frac{\nabla_{11}^2\psi}{b_{11}}+\frac{1-A}{h-h_t}+\frac Ah\notag\\
&\le&-A\sum b^{ii}-2B\sum b^{ii}(b_{ii}^2-2hb_{ii})+2B\sum h_k\Big( \frac{h_k}{h-h_t}+b^{kk}h_k-\nabla_k\psi\Big)\notag\\
&&- \frac{\nabla_{11}^2\psi}{b_{11}}+\frac{1-A}{h-h_t}+CA\notag\\
&\le&(2B|\nabla h|^2-A)\sum b^{ii}-2B\sum b_{ii} +\frac{1-A+2B|\nabla h|^2}{h-h_t}\\
&&- \frac{\nabla_{11}^2\psi}{b_{11}}-2B\sum h_k\nabla_k\psi+C(A+B).\notag
\eeqn

On the other hand, direct computation shows that
		\begin{equation}\label{psi 1}
			\nabla_k \psi=\frac{f_k}{f}-\frac{\nabla_k \widetilde{V}_{q-1}\left(\Omega_t, \overline{\nabla} h\right)}{\widetilde{V}_{q-1}\left(\Omega_t, \overline{\nabla} h\right)}+p\frac{h_k}{h},
		\end{equation}
and
		\begin{equation}\label{psi 11}
			\nabla^2_{11} \psi=\frac{f f_{11}-f_1^2}{f^2}-\frac{\nabla^2_{11} \widetilde{V}_{q-1}\left(\Omega_t, \overline{\nabla} h\right)}{\widetilde{V}_{q-1}\left(\Omega_t, \overline{\nabla} h\right)}+\frac{\left(\nabla_1 \widetilde{V}_{q-1}\left(\Omega_t, \overline{\nabla} h\right)\right)^2}{\left(\widetilde{V}_{q-1}\left(\Omega_t, \overline{\nabla} h\right)\right)^2}+p\frac{h h_{1 1}-h_1^2}{h^2}.
		\end{equation}
As a consequence of \eqref{psi 1} and \eqref{psi 11}, one infers that at $(x_0,t_0)$,
\beqn\label{eq psi0}	
&&- \frac{\nabla_{11}^2\psi}{b_{11}}-2B\sum h_k\nabla_k\psi\notag\\
&\le&\widetilde{V}^{-1}_{q-1}(\Omega_t, z_0) \Big( \frac{\nabla_{11} \widetilde{V}_{q-1}(\Omega_t, z_0)}{b_{11}} 
+2B\sum h_k\nabla_{k} \widetilde{V}_{q-1}(\Omega_t, z_0)\Big) +C(B+1),
\eeqn
where $z_0 = \overline\nabla h(x_0,t_0)$.
Since $q>3$, it follows by \cite[Lemma 5.3]{HHL reg} that at $(x_0,t_0)$,
		\begin{equation}\label{i V}
			\nabla_k\widetilde{V}_{q-1}\left(\Omega_t, z_0\right)=\frac{(q-1)(n-q+1)}{n} 
			b_{kk}|_{(x_0,t_0)}\int_{\Omega_t} \frac{(y-z_0) \cdot e_k}{|y-z_0|^{n+3-q}} \mathrm{d} y,
		\end{equation}
and
\beqn\label{11 V}
&&\frac{n}{(q-1)(n-q+1)} \nabla_{11}^2\widetilde V_{q-1}(\Omega_t,z_0)\notag\\
& &=b_{11 k}|_{(x_0,t_0)}\int_{\Omega_t} \frac{(y-z_0) \cdot e_k}{|y-z_0|^{n+3-q}} \mathrm{d} y 
-b_{11}|_{(x_0,t_0)}\int_{\Omega_t}\frac{(y-z_0) \cdot x_0}{ |y-z_0|^{n+3-q} } dy\notag \\
&&+b_{11}^2|_{(x_0,t_0)} \int_{\Omega_t}|y-z_0|^{q-n-3}\Big[(n+3-q) \frac{\left((y-z_0) \cdot e_1\right)^2}{|y-z_0|^2} -1\Big] \mathrm{d} y.
\eeqn
Plugging \eqref{i V} and \eqref{11 V} in \eqref{eq psi0}, we obtain
\beqs	
&&- \frac{\nabla_{11}^2\psi}{b_{11}}-2B\sum h_k\nabla_k\psi\notag\\
&\le&\frac{(q-1)(n-q+1)}{n\widetilde{V}_{q-1}(\Omega_t, z_0)}  
\Big(\frac{b_{11 k}}{b_{11}} +2Bh_k b_{kk} \Big)\int_{\Omega_t} \frac{(y-z_0) \cdot e_k}{|y-z_0|^{n+3-q}}dy \notag\\
&& +Cb_{11}+C(B+1).
\eeqs
Using \eqref{T3.1:cc1}, we further conclude that
\beqn\label{eq psi}	
- \frac{\nabla_{11}^2\psi}{b_{11}}-2B\sum h_k\nabla_k\psi \le Cb_{11}+C(B+1).
\eeqn
Inserting \eqref{eq psi} into \eqref{ptE} and choosing $A=2 B \max _{\mathbb{S}^{n-1} \times(0,+\infty)}|\nabla h|^2+1$,
\beqs
0\le -2B\sum b_{ii} +Cb_{11}+ C(A+B), \quad\text{at}\ (x_0,t_0).
\eeqs
Taking $B$ large, we conclude from the above inequality that $b_{11}(x_0,t_0)\le C$ as desired.
\end{proof}
	
	With the second derivative estimates \eqref{2th derivative}, the equation \eqref{flow2} are uniformly parabolic. 
	By \cite[Theorem 1.2]{HHL reg}, we have
	\begin{equation}\label{Vq C2estimate}
		\|\widetilde{V}_{q-1}\left(\Omega_t,\overline{\nabla}h\right)\|_{C^2\left(\mathbb{S}^{n-1}\right)} \leq C \quad \forall(x, t) \in \mathbb{S}^{n-1} \times[0, T).
	\end{equation}
Using the Krylov regularity theory \cite{Krylov} and a bootstrap argument, we obtain
	$$
	\|h(\cdot, t)\|_{C^{3, \alpha}\left(\mathbb{S}^{n-1}\right)} \leq C \quad \forall(x, t) \in \mathbb{S}^{n-1} \times[0, T),
	$$
for any given $\alpha \in(0,1)$, where the constant $C$ depends only on $\alpha, n, p, q$, $\min _{\mathbb{S}^{n-1}} f$, 
$\|f\|_{C^{1,1}\left(\mathbb{S}^{n-1}\right)}$, and the initial condition on $h(\cdot, 0)$. 
We hence conclude the long-time existence of solutions to the flow \eqref{flow1}.
	
\begin{theorem}\label{longtime}
Let $f \in C^{1,1}(\mathbb{S}^{n-1}) $ be a positive function and $T_{\max }$ be the maximal time such that 
$h(\cdot, t)$ is a positive, $C^{3, \alpha}$-smooth, and uniformly convex solution to \eqref{flow2} on $\left[0, T_{\max }\right)$. 
If $p<-n-q+1$ and $3< q<n+1$, and \eqref{C0C1} holds for all $t\in[0, T_{\max})$,
then $T_{\max }=\infty$.
\end{theorem}
	
	\section{Properties of the functional and the initial condition}
	In this section, we find a nice initial hypersurface $\mathcal{N}_0$ such that the flow \eqref{flow1} deforms $\mathcal{N}_0$ into a solution to \eqref{MAeq} after appropriate scaling.  We first demonstrate the monotonicity of the functional \eqref{functional} under the flow \eqref{flow1}.
	\begin{lemma}
Let $\mathcal{M}_t,$ $t\in [0,T),$ be a solution to the flow \eqref{flow1} in $\mathcal{K}_o.$ Then
		\begin{equation}\label{monotone}
			\frac{\mathrm{d}}{\mathrm{d} t} \mathcal{J}(\Omega_t)\geq 0,
		\end{equation}
		where $\Omega_t= \text{Cl} (\mathcal{M}_t)$ is the convex body enclosed by $\mathcal M_t$.
		 The equality holds if and only if the support function of $\Omega_t$ satisfies \eqref{MAeq} after appropriate scaling.
	\end{lemma}
	\begin{proof}
		Let $h(\cdot,t)$ be the corresponding support function of $\Omega_t.$ By \eqref{diff}, we have 
		\begin{eqnarray*}
			\frac{\mathrm{d}}{\mathrm{d} t} \mathcal{J}(K_t)&=& -\int_{\mathbb{S}^{n-1}} f h^{p-1} \partial_t h \mathrm{d} \sigma_{\mathbb{S}^{n-1}}+\int_{\mathbb{S}^{n-1}} \partial_t h \mathrm{d} F_q(K_t,\cdot)\\
			&=& \int_{\mathbb{S}^{n-1}}\left(\frac{2q \widetilde{V}_{q-1}\left(\Omega_t,\overline{\nabla}h \right) \det\left(\nabla ^2 h+ h I\right)}{\omega_n }- f h^{p-1}\right) \partial_t h \mathrm{d} \sigma_{\mathbb{S}^{n-1}}\\
			&=& \int_{\mathbb{S}^{n-1}}\left(\frac{1}{ \mathcal{K}}-\frac{\omega_{n}}{2q}\frac{f h^{p-1}}{\widetilde{V}_{q-1}\left(\Omega_t,\overline{\nabla}h \right)}\right)^2 h \mathcal{K} \frac{2q \widetilde{V}_{q-1}\left(\Omega_t,\overline{\nabla}h \right)}{\omega_n } \mathrm{d} \sigma_{\mathbb{S}^{n-1}}\geq 0.
		\end{eqnarray*}
		Moreover, we can see directly that the equality $\frac{\mathrm{d}}{\mathrm{d} t} \mathcal{J}(K_t)=0$ holds if and only if $\lambda h(\cdot, t)$ is a solution to \eqref{MAeq} with $\lambda =\left(\frac{2q}{\omega_n}\right)^{\frac{1}{n+q-1-p}}.$ 
	\end{proof}
	
	Next, we establish the following property: for any given positive constant $A$, if one of $e_{\Omega}$, $\Vol(\Omega)$, $[\Vol(\Omega)]^{-1}$, or $[\operatorname{dist}(0, \partial \Omega)]^{-1}$ is sufficiently large, then $\mathcal{J}(\Omega) > A.$
	\begin{lemma}\label{prop 2}
		Suppose that $p<-n-q+1$ and $1 / c_0 \leq f \leq c_0$ for some $c_0 \geq 1$. For any given constant $A$, there exists positive constants $\delta, v_0,v_1, \bar{e}>0$ depending only on $n, p, c_0,q$ and $A$ such that if $\Omega \in \mathcal{K}_o$ in $\mathbb{R}^{n}$ satisfies one of the four cases (i) $\operatorname{dist}(0, \partial \Omega) \in(0, \delta)$;
(ii) $\operatorname{Vol}(\Omega) \geq v_1$; (iii) $\operatorname{Vol}(\Omega) \leq v_0$; (iv) $e_{\Omega} \geq \bar{e}$; then
		$$
		\mathcal{J}(\Omega)>A.
		$$
	\end{lemma}
	\begin{proof}
We first discuss case (i).
Let $E$ be the minimum ellipsoid of $\Omega$. 
By a rotation of coordinates, we assume
		$$
		E-\xi_E=\left\{z \in \mathbb{R}^{n}: \sum_{i=1}^{n} \frac{z_i^2}{a_i^2} \leq 1\right\},
		$$
		where $\xi_E$ as the center of $E$.
Assume $0<a_1\leq a_2 \leq \cdots \leq a_n$.
		Denote 
		\[
		d=\operatorname{dist}(0, \partial \Omega).
		\] 
Let $x_0 \in \mathbb{S}^{n-1}$ be the point such that
		$$
		h(x_0)=\min _{\mathbb{S}^{n-1}} h=d,
		$$
		where $h$ is the support function of $\Omega$. 
		Choose $j_0$ such that $$x_0 \cdot \mathbf{e}_{j_0}=\max \left\{\left|x_0 \cdot \mathbf{e}_i\right|: 1 \leq i \leq n\right\}.
		$$
Then $x_0 \cdot \mathbf{e}_{j_0} \geq c_n$, for some constant $c_n>0$ depending only on $n$. 
By John's lemma, $\frac{1}{n}E\subset \Omega\subset E$. Hence $a_1\geq c_nd.$
This together with Lemma \ref{elliptic estimate} yields
		\begin{equation}\label{chordeq}
			\mathcal{J}(\Omega)>I_q(\Omega) \geq I_q(\frac{1}{n}E)\geq c_n\prod_{i=2}^n a_i a_1^q \geq c_n d^q \prod_{i\neq j_0}^n a_i\geq c_n d^q \prod_{i\neq j_0}^n h(\mathbf e_i).
		\end{equation}
		
On the other hand, by the proof of  \cite[Lemma 2.2]{GLW2022} (the part between estimate (2.6) to (2.9) there), we have
		\begin{equation}\label{chordeq2}
			\mathcal{J}(\Omega)\geq -\frac{1}{p}\int_{\mathbb{S}^{n-1}}fh^p \mathrm{d}\sigma_{\mathbb{S}^{n-1}} \geq \frac{c_n}{c_0 d^{-p-n+1}}\left[\prod_{i \neq j_0}^n h\left(\mathbf{e}_i\right)\right]^{-1}.
		\end{equation}
Note that the authors consider the problem on $\mathbb S^n$ in \cite{GLW2022}, 
while in this paper we consider the problem on $\mathbb S^{n-1}$.
Combining \eqref{chordeq} and \eqref{chordeq2}, we obtain
		$$
		[\mathcal{J}(\Omega)]^2 \geq \frac{c_n}{c_0 d^{-p-n-q+1}}.
		$$
		Since $-p-n-q+1>0$, we see that $\mathcal{J}(\Omega)>A$ if $d<\delta=\delta(A)$.
		
We next consider case (ii).
Assume that $d=\operatorname{dist}(0, \partial \Omega) \geq \delta$, otherwise we are done. 
This implies $a_1\geq c_n \delta$. By Lemma \ref{elliptic estimate}, we have
		\begin{eqnarray*}
			I_q(\Omega) &\geq& c_n \prod_{i=2}^n a_i a_1^q\\
			&\geq & c_{n,q}\operatorname{Vol}(\Omega) a_1^{q-1}\\
			&\geq & c_{n,q}\operatorname{Vol}(\Omega) \delta^{q-1}.
		\end{eqnarray*}
Therefore $I_q(\Omega)$ is as large as we want if $\operatorname{Vol}(\Omega)$ is sufficiently large. 
Hence $\mathcal{J}(\Omega)>A$ provided $v_1$ is large. 
		
For case (iii), we employ \cite[claim 8.1]{LXZY2022} and find that
		$$
		I_q(\Omega)\leq c(n,q)\operatorname{Vol}(\Omega)^{\frac{n+q-1}{n}}\leq 
		c(n,q)v_0^{\frac{n+q-1}{n}}.
		$$
Since $B_d$, the ball centred at the origin with radius $d=\operatorname{dist}(0, \partial \Omega)$, lies in $\Omega$, we have
		$$
		c(n,q)v_0^{\frac{n+q-1}{n}} \geq I_q(\Omega) \geq I_q\left(B_d\right)=I_q\left(B_1\right) d^{n+q-1} .
		$$
If $v_0$ is sufficiently small, then $d<\delta$ with $\delta$ being the constant in case (i). We are done.
		
We finally discuss case (iv). Assume $d=\operatorname{dist}(0, \partial \Omega) \geq \delta$. Otherwise we are done.
Then $B_\delta \subset \Omega \subset E$. It follows that
		\begin{equation}\label{a1}
			\delta \leq d \leq C_n a_i \text{ for }i\in\{1,\cdots,n\}.
		\end{equation}
Using \eqref{a1} and Lemma \ref{elliptic estimate},
		$$
		\mathcal{J}(\Omega)>I_q(\Omega) \geq c_n I_q(E)\geq c_n \prod_{i=2}^n a_i a_1^q \geq c_n e_{\Omega} \delta^{n+q-1}\geq c_n \bar{e}\delta^{n+q-1} .
		$$
As a result, $\mathcal{J}(\Omega)>A$ if $\bar{e}$ is sufficiently large. 
	\end{proof}
	\begin{remark}\label{rmk}
		By Lemma \ref{prop 2}, we know that if $\mathcal{J}(\Omega_t)\leq A,$ then there are some constants, independent of $t,$ such that
		\begin{equation}\label{c0}
			e_{\Omega_t}\leq \bar{e},\quad v_0\leq \operatorname{Vol}(\Omega_t)\leq v_1,\quad\text{ and }B_{\delta}\subset \Omega_t.
		\end{equation}
		Note that this implies the $C^0$-estimate of $h(\cdot,t)$. 
	The $C^1$-estimate follows by the convexity of $h$. Therefore, all we need is the estimate $\mathcal{J}(\Omega_t)\leq A$ 
	for some constant $A$ independent of $t,$ where $\partial\Omega_t$ is a solution to \eqref{flow1}.
\end{remark}
	
Let us introduce the modified flow as in \cite{GLW2022}. Fix the constant
	\begin{equation}\label{A0}
		A_0=3I_q(B_1)+3 n^{-p}\|f\|_{L^1\left(\mathbb S^{n-1}\right)}
	\end{equation}
	such that if the minimum ellipsoid of $\Omega$ is $B_1(0)$, namely $\frac{1}{n} B_1(0) \subset \Omega \subset B_1(0),$ then
	$$
	\mathcal{J}(\Omega) \leq \frac{1}{2} A_0.
	$$
For a closed, smooth and uniformly convex hypersurface $\mathcal{N}$ such that $\Omega_0=\mathrm{Cl}(\mathcal{N}) \in \mathcal{K}_o$, we define $\overline{\mathcal{M}}_{\mathcal{N}}(t)$ with initial data $\mathcal{N}$ as follows:
	\begin{enumerate}
		\item If $\mathcal{J}\left(\mathcal{M}_{\mathcal{N}}(t)\right)<A_0$ for all time $t \geq 0$, let $\overline{\mathcal{M}}_{\mathcal{N}}(t)=\mathcal{M}_{\mathcal{N}}(t)$ for all $t \geq 0$, where $\mathcal{M}_{\mathcal{N}}(t)$ is the solution to the flow \eqref{flow1}. 
		\item If $\mathcal{J}(\mathcal{N})<A_0$, and $\mathcal{J}\left(\mathcal{M}_{\mathcal{N}}(t)\right)$ reaches $A_0$ at the first time $t_0>0$, we define
		$$
		\overline{\mathcal{M}}_{\mathcal{N}}(t)=\left\{\begin{array}{lr}
			\mathcal{M}_{\mathcal{N}}(t), & \text { if } 0 \leq t<t_0, \\
			\mathcal{M}_{\mathcal{N}}\left(t_0\right), & \text { if } t \geq t_0 .
		\end{array}\right.
		$$
		\item If $\mathcal{J}(\mathcal{N}) \geq A_0$, we let $\overline{\mathcal{M}}_{\mathcal{N}}(t) \equiv \mathcal{N}$ for all $t \geq 0$. That is, the solution is stationary.
	\end{enumerate}
	We call $\overline{\mathcal{M}}_{\mathcal{N}}$ a modified flow of \eqref{flow1}. By Lemma \ref{prop 2}, for the given constant $A_0$ in \eqref{A0}, there exist sufficiently small constants $\delta$ and $\bar v<1,$ and a sufficiently large constant $\bar{e}$ such that we have the following properties:
	\begin{enumerate}
		\item [(a):]If one of the four cases $\dist(0,\mathcal{N})<\delta,$ $\Vol(\mathrm{Cl}(\mathcal{N}))<\omega_n\bar v,$ $\Vol(\mathrm{Cl}(\mathcal{N}))>\frac{\omega_n}{\bar v},$ $e_{\mathrm{Cl}(\mathcal{N})}>\bar{e}$ occurs, we have $\mathcal{J}(\mathcal{N}) > A_0$ and so $\overline{\mathcal{M}}_{\mathcal{N}}(t) \equiv \mathcal{N}$ for all $t$.
		\item [(b):]If $\mathrm{Cl}(\mathcal{N})$ is very close to $B_1(0)$ in Hausdorff distance, then $\mathcal{J}(\mathcal{N})<A_0$.
		\item [(c):]By the definition of the modified flow, $\mathcal{J}\left(\overline{\mathcal{M}}_{\mathcal{N}}(t)\right)<\max \left\{A_0, \mathcal{J}(\mathcal{N})\right\}$ for all $t$.
	\end{enumerate}
	Hence, if $\overline{\mathcal{M}}_{\mathcal{N}}(t)$ is not identical to $\overline{\mathcal{M}}_{\mathcal{N}}(0)=\mathcal{N}$ for all $t>0$, then
	\begin{equation}\label{key c0}
		e_{\overline{\mathcal{M}}_{\mathcal{N}}(t)} \leq \bar{e}, \quad \omega_n\bar v \leq \operatorname{Vol}\left(\overline{\mathcal{M}}_{\mathcal{N}}(t)\right) \leq \omega_n\bar v^{-1}, \text { and } B_{\delta}(0) \subset \mathrm{Cl}\left(\overline{\mathcal{M}}_{\mathcal{N}}(t)\right) \quad \forall t \geq 0 .
	\end{equation}
	Denote 
	\begin{equation}\label{AI}
		\mathcal{A}_I=\{E\in \overline{\mathcal{K}_o} \text{ is an ellipsoid in }\mathbb{R}^n: \omega_n\bar v\leq \Vol(E)\leq \omega_n\bar v^{-1} \text{ and }e_E\leq \bar{e}\},
	\end{equation}

Since the chord integral and the eccentricity are invariant under the translation,  
it is convenient to consider ellipsoids centered at the origin. 
For every ellipsoid $E\in \mathcal{K}_o$, there exists a unique affine transformation $A$
(equivalently a positive definite matrix) such that $E=AB_1$. 
This observation together with \cite[Lemma 3.4]{GLW2022} implies that $\mathcal{A}_I$ is homeomorphic to $\mathcal{E}_I\times B_1$,
where
\begin{equation}\label{EI2}
	\mathcal{E}_I=\{A\in M^{n\times n}| A \text{ is positive definite }, \bar{v}\leq \det A\leq \frac{1}{\bar{v}},e_{A}\leq \bar{e}\},
\end{equation}
and $e_A$ denotes the ratio between the maximum eigenvalue and minimum eigenvalue of $A.$ 
Note that the eigenvalues of the matrix $A$ are the principal radii of $E$ and so $e_A=e_E$. 
Let $\mathcal{P}$ be the boundary of $\mathcal{A}_I\simeq\mathcal{E}_I\times B_1$.


\begin{lemma}[Lemma 3.5 in \cite{GLW2022}]\label{Psi retraction}
There is a retraction $\Psi$ from $\mathcal{A}_I\backslash\{B_1\}$ to $\mathcal{P}$.
\end{lemma}

Since  $\mathcal{A}_I$ is  homeomorphic to $\mathcal{E}_I\times B_1,$  the above lemma implies that
there exists a retraction from $(\mathcal{E}_I\times B_1)\backslash\{(I, 0)\}$ to $\mathcal{P}$,
where $I$ is the identity matrix. For simplicity of notations, we still use $\Psi$ to denote this retraction.

Instead of calculating the homology of $\mathcal{P}$ as in \cite{GLW2022}, 
we next apply the Brouwer fixed theorem to deduce the following key conclusion.  
This simplified the argument in \cite{GLW2022}.

\begin{lemma}\label{initial for t}
For every $t>0$, there exists $\mathcal{N}=\mathcal{N}_t$ with $\operatorname{Cl}(\mathcal{N}) \in \mathcal{A}_I$, such that the minimum ellipsoid of $\overline{\mathcal{M}}_{\mathcal{N}}(t)$ is the unit ball $B_1(0)$ centered at the origin.
\end{lemma}
\begin{proof}
Suppose by contrary that there exists $t_0>0$, such that for any $\Omega\in \mathcal{A}_I$,
the minimum ellipsoid of $\Omega_{\mathcal{N}}(t_0):=\text{Cl}(\overline{\mathcal{M}}_{\mathcal{N}}(t_0))$,
denoted by $E_{\mathcal{N}}(t_0)$, is not the unit ball $B_1(0)$.
Here $\mathcal{N}=\partial \Omega$.
By this assumption, there is a continuous map $T:\mathcal{A}_I\rightarrow \mathcal{A}_I\backslash \{B_1\}$ given by
$$ \mathcal{A}_I \ni\Omega\mapsto E_{\mathcal{N}}(t_0)\in \mathcal{A}_I\backslash \{B_1\}.$$ 
By Lemma \ref{prop 2} and the construction of the modified flow,
$T=id$ if restricted to $\mathcal{P} = \partial \mathcal A_I$.
Since  $\mathcal{A}_I$ is  homeomorphic to $\mathcal{E}_I\times B_1$, 
it implies that there exists a continuous map  $$\tilde{T}: \mathcal{E}_I\times B_1\rightarrow  (\mathcal{E}_I\times B_1)\backslash\{(I, 0)\},$$ such that $\tilde{T}=id$ on $\mathcal{P}$. 
It follows that
\begin{equation}\label{l4.5:tt}
\tilde \Psi = \Psi\circ \tilde T: \quad \mathcal{E}_I \times B_1 \to \mathcal{P}
\end{equation}
is a retraction, where $\Psi$ is the map in Lemma \ref{Psi retraction}.

Denote 
$$\mathcal{D}:=\{A\in M^{n\times n}| A \text{ is positive definite }, \|A\|_{\infty}\leq L,e_{A}\leq \bar{e}\}, $$  
where the constant  $L$ is chosen large so that $\mathcal{E}_I\subset \mathcal{D}.$
We claim that $\mathcal{D}$ is convex. Namely, if $A, B\in \mathcal D,$ 
then
\[
\lambda A+(1-\lambda)B\in \mathcal{D}\quad \text{for any}\ \lambda\in [0, 1].
\]  
For this end, write $A=\{a_{ij}\}^{n}_{i,j=1}$ and $B=\{b_{st}\}^{n}_{s,t=1}$. 
Take $C=\{c_{ij}\}^{n}_{i,j=1},$ where
\[
c_{ij}=\lambda a_{ij}+(1-\lambda)b_{ij}.
\]
It is clear that $|c_{ij}|\leq \lambda \|A\|_{\infty}+(1-\lambda)\|B\|_{\infty}\leq L.$ Let $a_1\leq a_2\le \cdots \leq a_n$
(resp. $b_1\leq b_2\le \cdots\leq b_n$) be the eigenvalues of $A$ (resp. $B$). 
Then it is straightforward
to check that the ratio between the maximum and the minimum eigenvalues of 
$C$ is bounded by
\[\frac{\lambda a_n+(1-\lambda)b_n}{\lambda a_1+(1-\lambda)b_1}\leq \bar e.
\] 
Hence $C\in \mathcal{D}$. Therefore $D$ is a convex set.


Now, we construct a retraction $\Phi:\mathcal{D}\times B_1\rightarrow \mathcal{E}_I\times B_1$ as follows:
given any $(A, z)\in \mathcal{D}\times B_1,$
\begin{itemize}
	\item[]\textbf{Case I:} $\bar{v}\leq \det A\leq \frac{1}{\bar{v}}$.  Take $\Phi (A, z)=(A, z)$.
	\item[]\textbf{Case II:} $\det A< \bar{v}$. 
Since $d(t):=\det^{\frac{1}{n}}(tI+(1-t)A)$ is concave with respect to $t$, $d(0)<\bar{v}^{\frac{1}{n}}<1$ and $d(1)=1$, 
there exists a unique constant $t_A>0$ such that
\[
d(t_A)=(\bar{v})^{\frac{1}{n}},\ \text{and} \ d(t)<(\bar{v})^{\frac{1}{n}}\ \forall \ t\in [0,t_A).
\]
Since $\det (t_A I+(1-t_A)A)=\bar{v},$ and
the eccentricity of $t_A I+(1-t_A)A$ is less or equal than $e_A$, we find that
$t_A I+(1-t_A)A\in \mathcal{E}_I$. Define
\[
\Phi (A, z)=(t_A I+(1-t_A)A, z).
\]

\item[]\textbf{Case III:} $\det A> \frac{1}{\bar{v}}$.
Since $d(t):=\det^{-\frac{1}{n}}(tI+(1-t)A^{-1})$ is convex with respect to $t$,
$d(0)>(\bar{v})^{-\frac{1}{n}}>1$ and $d(1)=1$,  
there is a unique constant $t_A>0$ such that
\[
d(t_A)=(\bar{v})^{-\frac{1}{n}},\ \text{and}\ d(t)>(\bar{v})^{-\frac{1}{n}}\ \forall t\in [0,t_A).
\]  
Since $\det (t_A I+(1-t_A)A^{-1})^{-1}=\frac{1}{\bar{v}},$ and
the eccentricity of $ (t_A I+(1-t_A)A^{-1})^{-1}$ is less or equal than $e_A$,
one infers that $ (t_A I+(1-t_A)A^{-1})^{-1}\in \mathcal{E}_I$. Take
\[
\Phi (A, z)=( (t_A I+(1-t_A)A^{-1})^{-1}, z).
\]
\end{itemize}

Observe that
$(A, z)\in \mathcal{P}$ if and only if one of the following four scenarios happens:
$|z|=1,$ or $\det A=\bar v$ or $\det A=\frac{1}{\bar v}$ or $e_A=\bar e.$
Let
$g:\mathcal{P}\rightarrow \mathcal{P}$ be a continuous map defined by
$g(A,z)=(A^{-1}, -z).$ It is straightforward to check that $g$ has no fixed point

Consider the map $G:= i\circ g\circ\tilde\Psi\circ\Phi$,
 where $i: \mathcal{P}\rightarrow \mathcal{D}\times B_1$ is the inclusion
 and $\tilde \Psi$ is given by \eqref{l4.5:tt}. 
 Then, by the above constructions, $G:\mathcal{D}\times B_1\rightarrow \mathcal{D}\times B_1$ is a continuous and has no fixed point.
 This contradicts to the Brouwer fixed point theorem, as $\mathcal{D}\times B_1$ is convex.	
Hence, for every $t>0$, we can find $\mathcal{N}$ with $\text{Cl}(\mathcal{N})\in \mathcal{A}_I$, such that the minimum ellipsoid of $\overline{\mathcal{M}}_{\mathcal{N}}(t)$ is the unit ball $B_1(0)$.
\end{proof}

By Lemma \ref{initial for t}, for a sequence $t_k\rightarrow \infty$, we can take the initial datas $\mathcal{N}_k=\mathcal{N}_{t_k}$
such that the minimum ellipsoid of $\overline{\mathcal M}_{\mathcal N_k} (t_k)$ is $B_1(0)$. By Blaschke selection theorem, we have that $\mathcal{N}_k$ converges to a limit $\mathcal{N}_0$ such that 
$\text{Cl}(\mathcal{N}_0)\in \mathcal{A}_I$ up to a subsequence in Hausdorff distance.
By the choice of $A_0$, namely \eqref{A0}, 
$$
\mathcal{J}\left(\overline{\mathcal{M}}_{\mathcal{N}_k}(t_k)\right)\leq \frac{A_0}{2}.
$$
The construction of the modified flow and the monotonicity of \eqref{functional} then yield
\[\overline{\mathcal{M}}_{\mathcal{N}_k}(t)=\mathcal{M}_{\mathcal{N}_k}(t), \ \forall t\leq t_k.\]
Following the proof of \cite[Lemma 3.11]{GLW2022}, we obtain
\begin{equation}\label{func of initial}
\mathcal{J}\left({\mathcal{M}}_{\mathcal{N}_0}(t)\right)< A_0,\quad\forall t\ge0.
\end{equation}

\section{Proof of Theorem \ref{A} and Theorem \ref{B}}
In this section, we show the convergence of the flow \eqref{flow2} with initial data $\mathcal N_0$ found in the end of Section 4.
Let $\Omega_{\mathcal{N}_0}(t)=\operatorname{Cl}\left(\mathcal{M}_{\mathcal{N}_0}(t)\right)$ and $h(\cdot, t)$ be its support function. By \eqref{func of initial} and Lemma \ref{prop 2}, we have
$$
B_{\delta} \subset \Omega_{\mathcal{N}_0}(t), \quad \omega_n\bar{v} \leq \text{Vol}\left(\Omega_{\mathcal{N}_0}(t)\right) \leq \omega_n\bar{v}^{-1}\text { and } e_{\mathcal{M}_{\mathcal{N}_0}(t)} \leq \bar{e}, \text { for all } t \geq 0.
$$
Hence, there is a constant $C>0$ only depending on $n, p, q$ and the lower and upper bounds of $f$ such that
$$
\frac{1}{C} \leq h(x, t) \leq C, \quad \forall(x, t) \in \mathbb{S}^{n-1} \times[0, \infty).
$$
By the convexity of $h(x, t)$ we have $|\nabla h|(x, t) \leq C$ for all $(x, t) \in \mathbb{S}^{n-1} \times[0, \infty).$


By Theorem \ref{longtime}, $h(\cdot,t)$ is positive, $C^{3,\alpha}$-smooth, and uniformly convex for all time $t\ge0$.
We are now at a place to prove the Theorem \ref{A} and Theorem \ref{B}.

\begin{proof}[Proof of Theorem \ref{A}]
By \eqref{func of initial} and \eqref{monotone}, we have
$$
\int_0^{\infty} \mathcal{J}^{\prime}\left(\mathcal{M}_{\mathcal{N}_0}(t)\right) \mathrm{d} t \leq \limsup _{T \rightarrow \infty} \mathcal{J}\left(\mathcal{M}_{\mathcal{N}_0}(T)\right)-\mathcal{J}(\mathcal{N}_0) \leq A_0,
$$
which implies that there exists a sequence $t_i \rightarrow \infty$ such that
$$
\mathcal{J}^{\prime}\left(\mathcal{M}_{\mathcal{N}_0}(t_i)\right)=\int_{\mathbb{S}^{n-1}}\left(\frac{1}{ \mathcal{K}}-\frac{\omega_{n}}{2q}\frac{f h^{p-1}}{\widetilde{V}_{q-1}\left(\Omega_t,\overline{\nabla}h \right)}\right)^2 h \mathcal{K} \frac{2q \widetilde{V}_{q-1}\left(\Omega_t,\overline{\nabla}h \right)}{\omega_n }\vert_{t=t_i} \mathrm{d} \sigma_{\mathbb{S}^{n-1}} \rightarrow 0 .
$$
Passing to a subsequence, we obtain by Theorem \ref{C2} that $h\left(\cdot, t_i\right) \rightarrow h_{\infty}$ in 
$C^{3, \alpha}(\mathbb{S}^{n-1})$-topology and $\lambda h_\infty$ is a solution to \eqref{MAeq} with $\lambda =\left(\frac{2q}{\omega_n}\right)^{\frac{1}{n+q-1-p}}.$ 
\end{proof}

Theorem \ref{B} is a consequence of the combination of Theorem \ref{A} and an approximation argument. 
By the weak convergence of the $L_p$ chord measure, the proof is  the same to that of  \cite[Corollary 1.2]{GLW2022}. 

\textbf{Acknowledgment.} Research of the first author was supported by National Key R\&D program of China (No.2022YFA1005400),  National Science Fund for Distinguished Young Scholars (No. 12225111), NSFC No.12141105.

\bibliographystyle{amsplain}

\end{document}